\documentclass{amsart}
\begin{document}

% ----------------------------------------------------------------
\vfuzz2pt % Don't report over-full v-boxes if over-edge is small
\hfuzz2pt % Don't report over-full h-boxes if over-edge is small
%Theorems
\newtheorem{thm}{Theorem}[section]
\newtheorem{corollary}[thm]{Corollary}
\newtheorem{lemma}[thm]{Lemma}
\newtheorem{proposition}[thm]{Proposition}
\newtheorem{defn}[thm]{Definition}
\newtheorem{remark}[thm]{Remark}
\newtheorem{example}[thm]{Example}
\newtheorem{fact}[thm]{Fact}
%\numberwithin{equation}{section}
\
% MATH -----------------------------------------------------------
\newcommand{\norm}[1]{\left\Vert#1\right\Vert}
\newcommand{\abs}[1]{\left\vert#1\right\vert}
\newcommand{\set}[1]{\left\{#1\right\}}
\newcommand{\Real}{\mathbb R}
\newcommand{\eps}{\varepsilon}
\newcommand{\To}{\longrightarrow}
\newcommand{\BX}{\mathbf{B}(X)}
\newcommand{\A}{\mathcal{A}}
\newcommand{\onabla}{\overline{\nabla}}
\newcommand{\hnabla}{\hat{\nabla}}

% ----------------------------------------------------------------

\def\proof{\medskip Proof.\ }
\font\lasek=lasy10 \chardef\kwadrat="32 %kwadrat
\def\kwadracik{{\lasek\kwadrat}}
\def\koniec{\hfill\lower 2pt\hbox{\kwadracik}\medskip}

\newcommand*{\C}{\mathbf{C}}
\newcommand*{\R}{\mathbf{R}}
\newcommand*{\Z}{\mathbf {Z}}

\def\sb{f:M\longrightarrow \C ^n}
\def\det{\hbox{\rm det}\, }
\def\detc{\hbox{\rm det }_{\C}}
\def\i{\hbox{\rm i}}
\def\tr{\hbox{\rm tr}\, }
\def\rk{\hbox{\rm rk}\,}
\def\vol{\hbox{\rm vol}\,}
\def\Im {\hbox{\rm Im}\, }
\def\Re{\hbox{\rm Re}\, }
\def\interior{\hbox{\rm int}\, }
\def\e{\hbox{\rm e}}
\def\pu{\partial _u}
\def\pv{\partial _v}
\def\pui{\partial _{u_i}}
\def\puj{\partial _{u_j}}
\def\puk{\partial {u_k}}
\def\div{\hbox{\rm div}\,}
\def\Ric{\hbox{\rm Ric}\,}
\def\r#1{(\ref{#1})}
\def\ker{\hbox{\rm ker}\,}
\def\im{\hbox{\rm im}\, }
\def\I{\hbox{\rm I}\,}
\def\id{\hbox{\rm id}\,}
\def\exp{\hbox{{\rm exp}^{\tilde\nabla}}\.}
\def\cka{{\mathcal C}^{k,a}}
\def\ckplusja{{\mathcal C}^{k+1,a}}
\def\cja{{\mathcal C}^{1,a}}
\def\cda{{\mathcal C}^{2,a}}
\def\cta{{\mathcal C}^{3,a}}
\def\c0a{{\mathcal C}^{0,a}}
\def\f0{{\mathcal F}^{0}}
\def\fnj{{\mathcal F}^{n-1}}
\def\fn{{\mathcal F}^{n}}
\def\fnd{{\mathcal F}^{n-2}}
\def\Hn{{\mathcal H}^n}
\def\Hnj{{\mathcal H}^{n-1}}
\def\emb{\mathcal C^{\infty}_{emb}(M,N)}
\def\M{\mathcal M}
\def\Ef{\mathcal E _f}
\def\Eg{\mathcal E _g}
\def\Nf{\mathcal N _f}
\def\Ng{\mathcal N _g}
\def\Tf{\mathcal T _f}
\def\Tg{\mathcal T _g}
\def\diff{{\mathcal Diff}^{\infty}(M)}
\def\embM{\mathcal C^{\infty}_{emb}(M,M)}
\def\U1f{{\mathcal U}^1 _f}
\def\Uf{{\mathcal U} _f}
\def\Ug{{\mathcal U} _g}
\def\[f]{{\mathcal U}^1 _{[f]}}
\def\hnu{\hat\nu}
\def\gnu{\nu_g}
\def\sup{\hbox{\rm sup}}
\def\oRic{\hbox{\overline\Ric}}
\def\hRic{\hbox{widehat{\Ric}}}

% WERSAJA OSTATECZNA 06.04.21

\title{Some inequalities and applications of Simons' type formulas in Riemannian, affine and statistical
geometry}
%\thanks
\author{Barbara Opozda}



%\thanks{The research supported by the NCN grant UMO-2013/11/B/ST1/02889
%and a grant of the TU  Berlin} \maketitle

\address{Faculty of Mathematics and Computer Science UJ,
ul. \L ojasiewicza  6, 30-348 Cracow, Poland}

\email{barbara.opozda@im.uj.edu.pl} \maketitle

\begin{abstract}
A few formulas and theorems for statistical structures are proved.
They deal with various curvatures as well as with metric properties
of the cubic form or its covariant derivative. Some of them
generalize formulas and theorems known in the case of Lagrangian
submanifolds or affine hypersurfaces.
% One of the
% purposes of this paper is to find general,
% unified, and optimal settings for all these formulas and theorems.
\end{abstract}
\section{Introduction}

The name "statistical geometry" is relatively new, although this
geometry  has existed  for long and in various editions. From a
geometric viewpoint, a statistical structure is nothing but a
Codazzi pair. Starting from locally strongly convex hypersurfaces in
the Euclidean space, through locally strongly convex equiaffine
hypersurfces in the affine space $\R^{n+1}$ and Lagrangian
submanifolds in complex space forms to Hessian manifolds -- all
these examples  are statistical manifolds. Note that the structures
of the subclasses have very different properties and, therefore, the
intersections of the subclasses are small. For instance, if a
statistical structure can be realized on a Lagrangian submanifold
and on an affine hypersurface, then its Riemannian sectional
curvature, its $\nabla$-sectional curvature and its $K$-sectional
curvature are all constant.  Note also that the category of
statistical structures is much larger than the union of  all the
specific subclasses mentioned above. In particular, results proved
for affine hypersurfaces or for Lagrangian submanifolds usually are
rarely  generalizable to the general case of statistical structures.

 The aim of this paper is
to prove some global and local theorems for statistical manifolds in
the  general setting, with extensive references to affine
hypersurfaces or Lagrangian submanifolds.

By a statistical structure we mean a pair $(g,\nabla)$, where $g$ is
a Riemannian metric tensor field and $\nabla$ is a torsion-free
affine connection such that the cubic form $\nabla g$ is symmetric.
Let  $A=-\frac{1}{2}\nabla g$ and $\tau$ be the Czebyshev form
defined as $\tau(X)=\tr_g A(X,\cdot,\cdot)$. Set $K=A^\sharp$. The
Levi-Civita connection for $g$ will be denoted by $\hat\nabla$.
$\hat R$ and $\widehat\Ric$ will stand for the corresponding
curvature and Ricci tensors. A statistical structure is called
trace-free if $\tau=0$. If $\hat\nabla A$ is symmetric, then the
statistical structure is called conjugate symmetric. Basic
information  on statistical structures, affine hypersurfaces, and
Lagrangian submanifolds, needed in this paper, is contained in
Section 2.

 In Section 3, we prove some algebraic formulas for
$\Vert A\Vert$, $\Vert \tau\Vert$, $\tau$, and $K$. These formulas
yield inequalities for the Ricci tensors and scalar curvatures of
statistical structures.

The main part of this paper deals with the following Simons' formula
\begin{equation}\label{Simons'_wstep}
\frac{1}{2}\Delta(\Vert s\Vert^2)=g(\Delta  s,s) +\Vert\hat\nabla
s\Vert^2,
\end{equation}
where $s$ is a tensor field and $\Delta s$ is a specially defined
Laplacian on $s$. The first assumption which should be made for
computing the term $g(\Delta s,s)$ is the symmetry of $\hnabla s$.
Therefore, when we apply this formula to the cubic form $A$, we
shall assume that the statistical structure is conjugate symmetric.
Using various properties of statistical structures, this term can be
computed in various ways. For example, one can get the equality, see
(\ref{2dla A}),

\begin{equation}\label{2dla A-wstep}
\frac{1}{2}\Delta (\Vert A\Vert^2)=\Vert\hnabla A\Vert^2+
g(\hnabla^2 \tau ,A)+g(\hat R-R,\hat
R)+g(\widehat{\Ric},g(K_\cdot,K_\cdot)),
\end{equation}
where $R$ is the curvature tensor of $\nabla$. Computing this term
in another way  one gets the following inequalities valid for
conjugate symmetric trace-free statistical structures, see
(\ref{inequalitiesNSLI}),
\begin{equation}\label{inequalitiesNSLI-wstep}
\begin{array}{rcl}
&&\frac{1}{2}\Delta \Vert A\Vert^2\ge (n+1)H\Vert A\Vert^2
+\frac{n+1}{n(n-1)}\Vert A\Vert^4+\Vert\hat\nabla
A\Vert^2,\\
&& \frac{1}{2}\Delta \Vert A\Vert^2\le (n+1)H\Vert
A\Vert^2+\frac{3}{2}\Vert A\Vert^4+\Vert\hat\nabla A\Vert^2
\end{array}
\end{equation}
in the case where $R=HR_0$, $R_0(X,Y)Z=g(Y,Z)X-g(X,Z)Y$, and $\dim
M=n$.
% We feature a few formulas of these types.
In Section 5, we employ (\ref{inequalitiesNSLI-wstep}) for proving
some estimations  for the functions $\Vert A\Vert ^2$,
$\Vert\hat\nabla A\Vert^2$ in the case where $g$ is complete. A
model theorem here is the following result due to Calabi
\begin{thm} For a complete hyperbolic  affine Blaschke sphere,
the Ricci tensor of $g$ is negative semi-definite and, consequently,
\begin{equation}
\Vert A\Vert ^2\le -\rho,
\end{equation}
where $\rho$ is the affine scalar curvature.
\end{thm}
In the applications of Simons' formula (\ref{Simons'_wstep}) the
term $\Vert\hat\nabla s\Vert^2$ is usually just nonnegative. We
propose  new  estimations involving  also  $\inf\Vert\hnabla
A\Vert^2$ and $\sup\Vert\hnabla A\Vert^2$, see Theorems
\ref{aboutinfu}, \ref{ostatnie}, and \ref{najostatnie}.
 For proving
these theorems, we use Yau's maximum principle. For example, we
prove (see Theorem \ref{ostatnie}) that if for a trace-free
statistical structure with complete metric on an $n$-dimensional
manifold we have $R=HR_0$, where $H$ is a negative number, then

\begin{equation}\label{infN_4}
\inf\, \Vert\hat\nabla A\Vert^2\le\frac{n(n^2-1)H^2}{4}.
\end{equation}

In the last section, we prove a theorem on conjugate symmetric
statistical structures with nonnegative sectional curvature of the
metric $g$. The theorem generalizes a theorem known for minimal
Lagrangian submanifolds, see \cite{MRU}.

\bigskip

\section{Preliminaries}

In this section, we  introduce the notation  and collect basic
information on the geometry of statistical structures, affine
hypersurfaces, and Lagrangian submanifolds. All details for this
section can be found in \cite{BW4}, \cite{NS}, and \cite{CHBY1}.

In this paper we consider only torsion-free connections. If needed,
 manifolds are connected and orientable. If $g$ is a metric tensor
field (positive-definite) on a manifold $M$, then $\div$ will stand
for a divergence relative to the Levi-Civita connection $\hnabla$ of
$g$. In particular, if $s$ is a tensor field of type $(1,k)$, then
$\div s$ is a tensor field of type $(0, k)$ defined as
\begin{equation}
(\div s)(X_1,...,X_k)=\tr\{Y\to(\hnabla _Ys )(X_1,...,X_k)\}.
\end{equation}
We shall  use  the  following Laplacian relative to $g$ acting on
tensor fields (usually symmetric) of type $(0,k)$ for various values
of $k$. Namely, we set
\begin{equation}\label{laplacianfors}
(\Delta s)(X_1,...,X_k)=\tr_g(\hnabla^2s)(\cdot,\cdot,
X_1,...,X_k)=\sum_i(\hnabla^2s)(e_i,e_i, X_1,...,X_k),
\end{equation}
where $e_1,...,e_n$ is an orthonormal basis of $T_xM$ and $
X_1,...,X_n\in T_xM$, $x\in M$. Moreover, $\hnabla^2
s=\hnabla(\hnabla s)$, where $\hnabla s$ is a tensor field of type
$(0,k+1)$ given by the formula
$$\hnabla s(X,... )=(\hnabla_{X} s)(...).$$
In the case of differential forms, the above Laplacian is  not the
Hodge Laplacian (defined as $d\delta+\delta d$). The two Laplacians
are related  by the Weitzenb${\rm\ddot{o}}$ck formulas, see e.g.,
(\ref{Weitzenbockfor1form}).
% Only in the case of functions both Laplacians coincide.
%, see e.g.\cite{BW4} and references there.
In this paper, the codifferential $\delta$  of a differential
$k$-form $\omega$ is defined as
\begin{equation}
\delta \omega(X_1,...,X_{k-1})=+ \tr_g(\hnabla w)(\cdot,\cdot,
X_1,..., X_{k-1}).
\end{equation}
Because of these agreements, Simons' formula has the form as in
(\ref{Simons'_wstep}).
 The Laplacian defined in
(\ref{laplacianfors}) is usually denoted by $\hnabla^*\hnabla$, but
we shall use the symbol $\Delta $ for  a simpler notation.

%The symbol $\Delta$ was also used in the sense of
%(\ref{laplacianfors}), for instance, in \cite{NS}.

\subsection{Statistical structures}
By a statistical structure on a manifold $M$ we mean a pair
$(g,\nabla)$, where $g$ is a metric tensor and $\nabla$ is a
torsion-free connection on $M$ such that the cubic form $\nabla g$
is symmetric.  The difference tensor defined by
\begin{equation}\label{K_nabla}
K(X,Y)=K_XY=\nabla_XY-\hat\nabla_XY
\end{equation}
 is symmetric and symmetric relative to $g$. A statistical structure is called trivial
 if $\nabla=\hat\nabla$, that is, $K\equiv 0$.
  It is
clear that a statistical structure can be also defined as a pair
$(g,K)$, where $K$ is a symmetric and symmetric relative to $g$
$(1,2)$-tensor field. Namely, the connection $\nabla$ is defined by
(\ref{K_nabla}). Alternatively, instead of $K$ one can use the
symmetric cubic form $A(X,Y,Z)=g(K(X,Y)Z)$. The cubic forms $\nabla
g$ and $A$ are related as follows
\begin{equation}
\nabla g=-2A.
\end{equation}
A statistical structure is called trace-free if $E:=\tr_gK=0$.
Having a metric tensor $g$ and a connection $\nabla$ on a manifold
$M$, one defines a conjugate connection $\onabla$ by the formula
\begin{equation}
g(\nabla _XY,Z)+g(Y,\onabla _XZ)=Xg(Y,Z).
\end{equation}
A pair $(g,\nabla)$ is a statistical structure if and only if
$(g,\onabla)$ is. The pairs are also simultaneously trace-free
because if $K$ is the difference tensor for $(g,\nabla)$, then $-K$
is the difference tensor for $\onabla$.

The curvature tensors for $\nabla$, $\onabla$, $\hnabla$ will be
denoted by $R$, $\overline R$ and $\hat R$ respectively. The
corresponding Ricci tensors will be denoted by $\Ric$,
$\overline\Ric$ and $\widehat\Ric$. In general, the curvature tensor
$R$ does not satisfy  the equality $g(R(X,Y)Z,W)=-g(R(X,Y)W,Z)$. If
it does, we say that the  statistical structure is conjugate
symmetric. For a statistical structure, we always have
\begin{equation}\label{RioR}
g(R(X,Y)Z,W)=-g(\overline  R(X,Y)W,Z).
\end{equation}
 The following
conditions are equivalent:
\newline
{\rm 1)} $R=\overline R$,
\newline
{\rm 2)} $\hnabla K$ is symmetric (equiv. $\hat\nabla A$ is
symmetric),
\newline
{\rm 3)} $g(R(X,Y)Z,W)$ is skew-symmetric relative to $Z,W$,

Hence, each of the above three conditions characterizes a conjugate
symmetric statistical structure. The above equivalences follow from
the following well-known formula
\begin{equation}\label{from_Nomizu_Sasaki}
R(X,Y)=\hat R(X,Y) +(\hnabla_XK)_Y-(\hnabla_YK)_X+[K_X,K_Y].
\end{equation}
Writing the same equality for $\onabla$ and adding both equalities,
we get
\begin{equation}\label{R+oR}
R+\overline R =2\hat R +2[K,K] ,\end{equation} where $[K,K]$ is a
$(1,3)$-tensor field defined as $[K,K](X,Y)Z=[K_X,K_Y]Z$.
 In particular, if $R=\overline R$, then
\begin{equation}\label{z[K,K]dlahiperpowierzchni}
R=\hat R +[K,K].
\end{equation}

The $(1,3)$-tensor field $[K,K]$ is a curvature tensor, that is, it
satisfies the first Bianchi identity and has the same symmetries as
the Riemannian curvature tensor. Define its Ricci tensor $\Ric^K$ as
follows $\Ric^K(Y,Z)=\tr\{X\to[K,K](X,Y)Z\}$. In \cite{BW4}, the
following formula   was proved
\begin{equation}\label{Ric^K}
\Ric^K (Y,Z)=\tau(K(Y,Z))-g(K_Y,K_Z),
\end{equation}
where  $\tau$ is the Czebyshev $1$-form defined as $\tau(X)=\tr
K_X=g(E,X)$. From (\ref{from_Nomizu_Sasaki}) we receive, see
\cite{BW4},

\begin{equation}\label{1zBW4}
\Ric=\widehat \Ric+\div K-\hat\nabla \tau +\Ric^K.
\end{equation}

The $1$-form $\tau$ is closed (i.e., $\hat\nabla\tau$ is symmetric)
if and only if the Ricci tensor $\Ric$ is symmetric. For a conjugate
symmetric statistical structure, the Ricci tensor is  symmetric.
More generally, if for a statistical structure
 $\Ric=\overline\Ric$, then
$d\tau=0$. Indeed, by writing (\ref{1zBW4}) for $\overline\nabla$
and comparing with  (\ref{1zBW4}), we see that $\Ric=\overline\Ric$
if and only if $\div K=\hat\nabla\tau$. Therefore, if
$\Ric=\overline\Ric$, then $d\tau=0$ (because $\div K$ is
symmetric). Using (\ref{1zBW4}) and the analogous formula for the
connection $\overline\nabla$ one also gets
\begin{equation}\label{Ric+oRic}
\Ric +\overline{\Ric} =2\widehat{\Ric}+2\tau\circ
K-2g(K_\cdot,K_\cdot).
\end{equation}
In particular, if $(g, \nabla )$ is trace-free, then
\begin{equation}\label{Riccihat_ge}
2\widehat{\Ric}\ge \Ric+\overline{\Ric}.\end{equation}

Recall now the scalar curvatures for statistical structures. We have
the Riemannian scalar curvature $\hat\rho$ for $g$ and the scalar
curvature $\rho$ for $\nabla$: $\rho=\tr_g\Ric$. It turns out (e.g.,
by (\ref{RioR})) that if we define the analogous scalar curvature
for $\onabla$, then it is equal to $\rho$. We also have the scalar
curvature $\rho^K=\tr_g \Ric^K$. From the above formulas, one easily
gets

\begin{equation}\label{rho^K}
\rho^K=\Vert E\Vert^2-\Vert K\Vert^2,
\end{equation}
\begin{equation}\label{theorema-egregium}
\hat\rho=\rho+\Vert K\Vert^2-\Vert E\Vert^2.
\end{equation}

Recall now the sectional curvatures in statistical geometry. Of
course, we have the usual Riemannian sectional curvature of the
metric tensor $g$. We shall denote it by $\hat k(\pi)$ if $\pi$ is a
vector plane in a tangent space. In \cite{BW4} (see also \cite{F})
another sectional curvature was introduced. It was called the
sectional $\nabla$-curvature. Namely, it was observed that the the
tensor field $R+\overline R$ has all algebraic properties needed to
produce the sectional curvature. However, in general, this sectional
curvature does not have the same properties as the Riemannian
sectional curvature. For instance, Schur's lemma does not hold, in
general. However, it holds for conjugate symmetric statistical
structures. The sectional $\nabla$-curvature is defined as follows
\begin{equation}
k(\pi)=\frac{1}{2}g((R+\overline R)(e_1,e_2)e_2,e_1),
\end{equation}
where $e_1,e_2$ is an orthonormal basis of $\pi$. Another good
curvature tensor, which exists in statistical geometry is $[K,K]$.
It again has  all algebraic properties to define a sectional
curvature. This sectional curvature was introduced and studied in
\cite{seccur}.  It was called the sectional $K$-curvature. If $\pi$
is a vector plane in a tangent space $T_xM$, then the sectional
$K$-curvature is equal to $g([K,K](e_1,e_2)e_2,e_1)$ for any
orthonormal basis $e_1,e_2$ of $\pi$. Similarly to the case of the
sectional $\nabla$-curvature Schur's lemma holds for conjugate
symmetric statistical structures.

A statistical structure is called Hessian if $\nabla$ is flat, that
is, if  $R=0$. As it will be noticed below, all  Hessian  structures
are locally realizable on parabolic equiaffine spheres. Since
$\Ric=0$ for a Hessian structure and $\nabla$ is torsion-free, we
know that the Koszul form $\beta=\nabla \tau$ is symmetric. A
Hessian structure is called  Einstein-Hessian if $\nabla\tau=\lambda
g$, see \cite{shima}. For a Hessian structure we have
\begin{equation}
\beta= \hat\nabla\tau -\tau\circ K,\ \ \ \ \tr _g\beta= \delta\tau
-\Vert \tau\Vert^2.
\end{equation}
By (\ref{Ric+oRic}), we also have
\begin{equation}
\widehat \Ric=g(K_\cdot,K_\cdot)-\tau\circ K.
\end{equation}

In geometry, the richest  sources of statistical manifolds seem to
be the theory of affine hypersurfaces in $\R^{n+1}$ and that of
Lagrangian submanifolds in complex space forms. We shall shortly
recall the basic facts from these theories.
\medskip

\subsection{Equiaffine locally strongly convex hipersurfaces in
$\R^{n+1}$}
 For the theory of affine hypersurfaces, we refer to \cite{NS}. Let
$f:M\to \R^{n+1}$ be an immersed locally strongly convex
hypersurface of the affine space $\R^{n+1}$. Denote the standard
flat affine connection (the operator of   ordinary differentiation)
on $\R^{n+1}$ by $D$. Assume that the hypersurface is equipped with
a transversal vector field $\xi$ such that $D_X\xi$ is tangent to
the hypersurface for every $X\in T_xM$, $x\in M$. Such a tranversal
vector field is called equiaffine and  a hypersurface endowed with
such a tranversal vector field is called an equiaffine hypersurface.
All hypersurfaces considered in this paper will be equiaffine. By
the following formulas of Gauss and Weingarten, one gets the induced
connection $\nabla$, the second fundamental form $g$, and the shape
operator $S$
\begin{equation}\label{Gauss}
\begin{array}{rcl}
&&D_Xf_*Y=f_*(\nabla_XY) +g(X,Y)\xi,\\
&&D_X\xi=-f_*(SX).
\end{array}
\end{equation}
A hypersurface is locally strongly convex if and only if $g$ is
definite. If an equiaffine hypersurface is locally strongly convex,
the sign of $\xi$ is chosen in such a way that $g$ is positive
definite. It turns out that for an equiaffine hypersurface the cubic
form $\nabla g$ is symmetric. The symmetry of $\nabla g$ is the
so-called first Codazzi equation for an equiaffine hypersurface. It
means that the induced structure $(g,\nabla)$ is a statistical
structure. We also have the following Gauss equation
\begin{equation}\label{Gauss_equation_for R}
R(X,Y)Z=g(Y,Z)\mathcal SX-g(X,Z)\mathcal SY,
\end{equation}
The so-called Ricci  equation for equiaffine immersions states that
the Ricci tensor $\Ric$ is symmetric. The symmetry of the Ricci
tensor is the first obstruction for a statistical structure to be
realizable as the induced structure on a hypersurface in $\R^{n+1}$.
The second obstruction is that the dual connection $\onabla$ must be
projectively flat. Namely, an important theorem on the realizability
of statistical structures on equiaffine hypersurfaces is the
following fundamental theorem proved in \cite{DNV}

\begin{thm}\label{fundamental}
Let $(g,\nabla)$ be a statistical structure  on a simply connected
manifold $M$ and satisfy the following conditions:
\newline
1) the Ricci tensor of $\nabla$   is symmetric,
\newline
2)  the dual connection  $\onabla$ is projectively flat.
\newline
Then there is a locally strongly convex immersion $f:M\to \R ^{n+1}$
and its equiaffine transversal vector field $\xi$ such that $\nabla$
is the induced connection and $g$ is the second fundamental form for
the equiaffine hypersurface $(f,\xi)$.
\newline
 \end{thm}

From now on, we shall always assume that a hypersurface is locally
strongly convex without mentioning this each time.

 An equiaffine hypersurface is called an equiaffine sphere
if the shape operator is a multiple of the identity, i.e., $S=H
\id$, where $\id$ is the identity $(1,1)$-tensor field on $M$. In
such a case, $H$ must be constant on  connected $M$. A hypersphere
is called elliptic if $H>0$, hyperbolic if $H<0$, and parabolic if
$H=0$. For an equiaffine sphere we have the equality $R=HR_0$, where
$R_0$ is the curvature tensor defined by $
R_0(X,Y)Z=g(Y,Z)X-g(X,Z)Y. $ It is now clear that for an equaiaffine
sphere its statistical structure is conjugate symmetric. The
converse is also true, that is, if the statistical structure on an
equiaffine  hypersurface is conjugate symmetric, then the
hypersurface must be a sphere (Lemma 12.5 in \cite{BW4}).

 It is also known that if we have a statistical
structure (on an $n$-dimensional connected manifold $M$) whose
curvature tensor $R$ satisfies the equality $R=HR_0$, where $H$ is
possibly a function, and  $n>2$, then $H$ must be a constant
(Proposition 12.7 in \cite{BW4}). If $R=HR_0$ then $g(R(X,Y)Z,W)$ is
skew-symmetric for $Z,W$, hence $\overline R=R=HR_0$. If moreover,
$H$ is constant $\onabla$ is clearly projectively flat. A
statistical structure for which $R=HR_0$, where $H$ is a constant,
was called in \cite{K} a statistical structure of constant
curvature. By Theorem \ref{fundamental}, we now have that
statistical structures of constant curvature are (from a local
viewpoint) exactly the induced structures on equiaffine spheres. Let
us rewrite (\ref{z[K,K]dlahiperpowierzchni}) in the case where
$R=HR_0$:
\begin{equation}\label{z[K,K]dlasfery}
HR_0=\hat R +[K,K].
\end{equation}
Among the equiaffine hypersuraces, the  best known and historically
 first discovered are those which additionally satisfy the so-called
  apolarity condition. This condition is equivalent to the
trace-freeness of the induced statistical structure (in the
terminology we are using in this paper). The importance of this
class of equiaffine hypersurfaces follows from the following
classical fact: {\it For a locally strongly convex hypersurface,
there is a unique (up to a constant) equiaffine transversal vector
field, for which the induced statistical structure is trace-free.}
An equiaffine hypersurface whose induced statistical structure is
trace-free is called a Blaschke hypersurface. An equiaffine sphere
whose induced structure is trace-free is called a Blaschke affine
sphere. In contrast to Riemannian geometry, the category of affine
spheres is very rich and not well recognized by now. Trivial
Blaschke hypersurfaces are quadrics. In particular, a parabolic
affine sphere with vanishing cubic form, that is, whose induced
statistical structure is trivial, must be  a part of an elliptic
paraboloid.

\medskip
\subsection{Lagrangian submanifolds}
We shall now briefly recall some facts and introduce some notations
for Lagrangian submanifolds. For this part we refer to \cite{CHBY1}.
Let $\tilde M(4c)$ be a complex space form of holomorphic sectional
curvature $4c$ and let $M$ be its Lagrangian submanifold, that is,
$M$ is a real submanifold of $\tilde M$, $\dim M=\dim_\C \tilde M$
and the bundle $JTM$ is orthogonal to $TM$, where $J$ is the complex
structure on $\tilde M$. The induced metric tensor field on $M$ will
be denoted by $g$. If $\sigma$ is the second fundamental tensor,
then we set $K=J\sigma$. The tensor field $K$ is symmetric and
symmetric relative to $g$. Hence, the pair $(g, K)$ is a statistical
structure. We shall call it the induced statistical structure on a
Lagrangian submanifold. For this structure we shall use all the
technique we use for statistical structure. The Codazzi equation for
a Lagrangian submanifold says that the $(1,3)$-tensor field $\hnabla
K$ is symmetric. Hence, the statistical structure on a Lagrangian
submanifold in a complex space form is conjugate symmetric. We have
the following Gauss equation for  a Lagrangian submanifold in
$\tilde M(4c)$
\begin{equation}\label{gausslagrangian}
cR_0=\hat R-[K,K].
\end{equation}
This equality is different than the  equality (\ref{z[K,K]dlasfery})
holding for equiaffine spheres. By comparing the two equalities, we
see that a statistical structure on  a Lagrangian submanifold can be
also realized (locally) on an equiaffine hypersurface  only when the
Riemannian sectional curvature for $g$ and the sectional
$K$-curvature are constant.
%In fact, one can easily see (using
%suitable fundamental theorems) that the converse is also true. More
%precisely, if for a given statistical structure all sectional
%curvatures, that is, the sectional curvature for $g$, the sectional
%$\nabla$-curvature, the sectional $K$-curvature are constant, then
%the structure can be locally realizable on an equiaffine sphere and
%on a Lagrangian submanifold.

Instead of the mean curvature vector of a Lagrangian submanifold,
one can consider the tangent vector field $E=\tr _gK$. In
particular, a Lagrangian submanifold of $\tilde M(4c)$ is minimal if
the induced statistical structure is trace-free. The mean curvature
tensor $\mu$ is equal to $-\frac{JE}{n}$. Hence it is parallel if
and only if $\hnabla E=0$. By (\ref{gausslagrangian}), we have
\begin{equation}
\widehat\Ric=c(n-1)g-g(K_\cdot,K_\cdot)+\tau\circ K
\end{equation}
and
\begin{equation}\label{theorema_egregium_Lagrangian}
\hat\rho =cn(n-1)+\Vert E\Vert^2-\Vert A\Vert ^2=cn(n-1)+n^2\Vert
\mu\Vert^2-\Vert A\Vert ^2.
\end{equation}

\bigskip

\section{Algebraic and curvature inequalities}

The notation $\sum _{i\ne j}$ stands for the sum over indices $i,j$
such that $i\ne j$, whereas $\sum_{j;j\ne i}$ means that the sum is
for the index $j$ such that $j\ne i$ for a fixed $i$.
\begin{thm}
For a statistical structure we have
\begin{equation}\label{nierownosc1}
(\tau\circ K)(U,U)-g(K_{U}, K_{U})\le \frac{1}{4}\Vert
\tau\Vert^2g(U,U)
\end{equation}
for every vector $U$. If $U\ne 0$, the equality holds if and only if
$\tau=0$ and $K_U=0$.

%For a unit vector $U$ the equality holds if and only if
%$K_UV=g(U,V)\lambda U$ and $\tau (V)=2\lambda g(U,V)$ for some
%number $\lambda$. If the equality holds for every $U$, then $K=0$.

If $g(K(U,U),U)=0$, then
\begin{equation}\label{nierownosc1a}
(\tau\circ K)(U,U)-g(K_{U}, K_{U})\le \frac{1}{8}\Vert
\tau\Vert^2g(U,U).
\end{equation}
If $U$ is unit, the equality holds if and only if $g(K(U,V), W)=0$
for $V,W$ perpendicular to $U$, $E$ is perpendicular to $U$ and
$E=4K(U,U)$.

\end{thm}
\proof We shall  prove  (\ref{nierownosc1}) for an arbitrary unit
vector $U\in T_xM$. We extend a given vector $U$ to an orthonormal
basis $e_1=U$, $e_2,..., e_n$ of $T_xM$. We set
$A_{ijk}=A(e_i,e_j,e_k)$ and $a_{ir}=A_{iir}$. We want to show
\begin{equation}\label{nierownosc3}
0\le
\frac{1}{4}\sum_{ijr}A_{iir}A_{jjr}+\sum_{ir}A_{1ir}^2-\sum_{ir}A_{iir}A_{11r}.
\end{equation}
The last formula can be rewritten as follows
\begin{equation}
0\le\frac{1}{4}\sum_{ijr}a_{ir}a_{jr} + a_{11}^2
+2\sum_{r=2}^na_{1r}^2+\sum_{i\ne1 r\ne 1}
A_{1ir}^2-\sum_{ir}a_{ir}a_{1r}.
\end{equation}
Since $\sum_{i\ne1 r\ne 1} A_{1ir}^2\ge 0$, it is sufficient to
prove

\begin{equation}\label{nierownosc4}
0\le a_{11}^2+2\sum_{r=2}^n
a_{1r}^2+\frac{1}{4}\sum_{ijr}a_{ir}a_{jr}-\sum_{ir}a_{ir}a_{1r}.
\end{equation}
The last inequality is equivalent to
\begin{equation}\label{nierownosc5}
0\le  a_{11}^2+2\sum_{r=2}^n a_{1r}^2
+\frac{1}{4}\sum_{rj;j\ne1}a_{1r}a_{jr}+\frac{1}{4}\sum_{r=1}^n\left
(\sum_{j=2}^n a_{jr}\right )^2 -\frac{3}{4}\sum_{ir}a_{ir}a_{1r}.
\end{equation}
For a fixed $r$, set $c_r=\sum_{j=2}^na_{jr}$.  We have
\begin{equation}\label{nierownosc5a}
0\le
\frac{1}{4}(a_{1r}-c_r)^2=a_{1r}^2+\frac{1}{4}a_{1r}c_r+\frac{1}{4}c_r^2-\frac{3}{4}a_{1r}(c_r+a_{1r}).
\end{equation}
 By comparing the last inequality with (\ref{nierownosc5}),
  we obtain
(\ref{nierownosc4}). Indeed, for $r=1$ we get in
(\ref{nierownosc5a}) $$ 0\le
a_{11}^2+\frac{1}{4}a_{11}c_1+\frac{1}{4}c_1^2-\frac{3}{4}a_{11}(c_1+a_{11}),
$$
hence this term
%$a_{11}^2+\frac{1}{4}a_{11}c_1+\frac{1}{4}c_1^2-\frac{3}{4}a_{11}(a_{11}+c_1)$
appearing also in (\ref{nierownosc5}) is nonnegative and equal to
$0$ if and only if $a_{11}=c_1$. Consider now $r>1$. For such  $r$
we have in (\ref{nierownosc5}) the term
\begin{equation}
2a_{1r}^2+\frac{1}{4}a_{1r}c_r+\frac{1}{4}c_r^2-\frac{3}{4}a_{1r}(c_r+a_{1r}).
\end{equation}
By comparing this with (\ref{nierownosc5a}), we see that this term
is nonnegative and equal to $0$ if and only if $a_{1r}=0$ and
$a_{1r}=c_r$. The inequality (\ref{nierownosc1}) has been proved and
moreover, one sees see that the equality holds if and only if
$A_{11r}=0$, $A_{22r}+...+A_{nnr}=0$ for every $r>1$,
$A_{111}=A_{221}+...+A_{nn1}$ and $A_{1ir}=0$ $\forall i\ne 1, r\ne
1$. In particular, $A_{1ii}=0$ for $i>1$. These conditions imply
that $K_{e_1}=0$ and $\tau=0$.

\medskip

Assume now that $g(K(U,U),U)=0$, that is, $a_{11}=0$. In
(\ref{nierownosc3}),  we replace $\frac{1}{4}$ by $\frac{1}{8}$ and
we obtain the following formula analogous to (\ref{nierownosc5})
\begin{equation}\label{nierownosc5b}
0\le  2\sum_{r=2}^n a_{1r}^2
+\frac{1}{8}\sum_{jr;j\ne1}a_{1r}a_{jr}+\frac{1}{8}\sum_{r=1}^n\left
(\sum_{j=2}^n a_{jr}\right )^2 -\frac{7}{8}\sum_{ir}a_{ir}a_{1r}.
\end{equation}
Again we fix $r$. The only possibly non-vanishing term of the right
hand side in (\ref{nierownosc5b}) with $r=1$ is $\frac{1}{8}c_1^2$.
For $r>1$ we write
\begin{equation}\label{nierownosc5c}
0\le
\frac{1}{8}(3a_{1r}-c_r)^2=2a_{1r}^2+\frac{1}{8}a_{1r}c_r+\frac{1}{8}c_r^2-\frac{7}{8}a_{1r}(c_r+a_{1r}).
\end{equation}
Comparing (\ref{nierownosc5c}) with (\ref{nierownosc5b}), we get the
desired inequality.

If we have the equality in (\ref{nierownosc1a}), then $c_1=0$. Since
$a_{11}=0$, we have $g(E,e_1)=0$.  By (\ref{nierownosc5c}), we have
$3A_{11r}=A_{22r}+...+A_{nnr}$. It follows that $E=4K(e_1,e_1)$. As
in the first part of the proof, we have that $A_{1ij}=0$ for
$i,j>1$. If $E=0$, then $K_{e_1}=0$.  Assume that $E\ne0$ and
$e_2=\frac{E}{\Vert E\Vert}$. If $K(e_1,e_1)=\lambda e_2$, then
$g(K_{e_1},K_{e_1})=2\lambda ^2$, $\tau(K(e_1,e_1))=4\lambda^2$ and
$\Vert\tau\Vert^2=16\lambda^2$.\koniec

As an immediate consequence of the inequality (\ref{nierownosc1})
and the formula (\ref{Ric+oRic}), we get a generalization of
(\ref{Riccihat_ge})
\begin{thm}
For any statistical structure we have
\begin{equation}
2\widehat\Ric\ge \Ric+\overline\Ric- \frac{1}{2}\Vert\tau\Vert^2g.
\end{equation}
For a statistical structure on an equiaffine sphere with $R=HR_0$,
we have
\begin{equation}
\widehat\Ric\ge \left (H(n-1)- \frac{1}{4}\Vert\tau\Vert^2\right)g.
\end{equation}
For a  Lagrangian submanifold of a complex space form $\tilde
M(4c)$, we have
\begin{equation}
\widehat\Ric\le \left(c(n-1)+\frac{n^2}{4}\Vert \mu\Vert ^2\right)g,
\end{equation}
where $\mu$ is the mean curvature vector of the submanifold.

\end{thm}
The last inequality was first proved in \cite{ChBY}.

\medskip
\medskip

\begin{thm} For any statistical structure, we have
\begin{equation}\label{nierownosc2}
\frac{n+2}{3}\Vert A\Vert^2-\Vert E\Vert^2\ge 0.
\end{equation}

\end{thm}
\proof We use the same notation as in the proof of the above
theorem.
 To prove (\ref{nierownosc2}), we first compute $\Vert A\Vert^2-\Vert
 E\Vert^2$. Let $e_1,...,e_n$ be any orthonormal basis of $T_xM$, $x\in M$. We have
\begin{equation}
\begin{array}{rcl}
&&\Vert E\Vert^2=\sum_i\left(\sum_j
A_{jji}\right)^2\\
&&\ \ \ \ =\sum_{i}A_{iii}^2+2\sum_{i\ne j}A_{iii}A_{jji}+\sum_{j\ne
i}A_{jji}^2+2\sum_{i\ne j, l\ne i, l< j}A_{jji}A_{lli}\\
&&\ \ \ \ \ \ \ =\sum_{i}a_{ii}^2+2\sum_{i\ne
j}a_{ii}a_{ji}+\sum_{j\ne i}a_{ji}^2+2\sum_{i\ne j, l\ne i, l<
j}a_{ji}a_{li}.
\end{array}
\end{equation}
\begin{equation}\label{a2}
\Vert A\Vert ^2= \sum_{ijl}A_{ijl}^2=\sum_{i}A_{iii}^2 +3\sum_{i\ne
j}A_{jji}^2 +\epsilon=\sum_{i}a_{ii}^2 +3\sum_{i\ne j}a_{ji}^2
+\epsilon,
\end{equation}
where $\epsilon=\sum_{i\ne j, i\ne r\j\ne i}A_{ijr}^2$. We now have
\begin{equation}\label{a2-b2}
\Vert A\Vert^2-\Vert
 E\Vert^2=2\sum_{i\ne j}a_{ji}^2-2\sum_{i\ne j, i\ne l,
 l<j}a_{li}a_{ji}-2\sum_{i\ne j}a_{ii}a_{ji}+\epsilon.
\end{equation}
For any real numbers $b_1,..., b_k$, the following formula holds
\begin{equation}
\sum_{l<j\le k}(b_l-b_j)^2=(k-1)\sum _jb_j^2-2\sum_{l<j\le k}b_lb_j.
\end{equation}
Hence, for a fixed $i$ we have
\begin{equation}\label{b}
(n-2)\sum_{j; j\ne i}a_{ji}^2-2\sum_{jl;j\ne i, l\ne i,
l<j}a_{li}a_{ji}\ge 0.
\end{equation}
Denote the left hand side of the last inequality by $c_i$ and set
$c=\sum_i c_i$. By (\ref{a2}), (\ref{a2-b2}), and (\ref{b}), we now
obtain

\begin{equation}\label{nierownosc6}
\frac{n-4}{3}\Vert A\Vert ^2+(\Vert A\Vert^2-\Vert E\Vert ^2)=
\frac{n-4}{3}\sum_i a_{ii}^2-2\sum_{j\ne
i}a_{ii}a_{ji}+c+\frac{n-1}{3}\epsilon.
\end{equation}
 For a fixed $i$ and  a positive number $\mu$, we have
\begin{equation}\label{mu}
\begin{array}{rcl}
&&0\le \sum _{j;j\ne i}\left( \mu
a_{ii}-\frac{1}{\mu}a_{ji}\right)^2=\sum _{j;j\ne i}\left[(\mu
a_{ii})^2-2(\mu a_{ii})\left(\frac{1}{\mu}
a_{ji}\right)+\left(\frac{1}{\mu}a_{ji}\right)^2\right]\\
&&\ \ \ \ \ \ \ \ = (n-1)\mu^2 a_{ii}^2-2\sum_{j, j\ne i}
a_{ii}a_{ji}+\frac{1}{\mu^2}\sum_{j; j\ne i}a_{ji}^2.
\end{array}
\end{equation}
We now add to the left hand side of (\ref{nierownosc6})  $\Vert
A\Vert^2$ and we set $\mu^2=\frac{1}{3}$. By (\ref{nierownosc6}),
(\ref{a2}), and (\ref{mu}), we obtain
\begin{equation}
\begin{array}{rcl}
&&\frac{n+2}{3}\Vert A\Vert^2 -\Vert E\Vert ^2=\Vert A\Vert^2+
\frac{n-4}{3}\Vert A\Vert^2 +(\Vert A\Vert^2-\Vert E\Vert^2)\\
&&\ \ \ \ \ \ =\frac{n-1}{3}\sum_ia_{ii}^2-2 \sum_{i\ne j}
a_{ii}a_{ji}+3\sum_{i\ne j} a_{ji}^2+\frac{n+2}{3}\epsilon+c\ge 0.
\end{array}
\end{equation}

Observe additionally that in the last line the equality holds if and
only if for every $i$ and $j$ such that $i\ne j$ we have
$A_{iii}=3A_{jji}$ and $A_{ijr}=0$ for $i\ne j, i\ne r, j\ne r$.

\koniec

 \begin{example}
 {\rm Consider  a statistical structure $(g,K)$ on $\R^2$, where $g$
 is the
 standard metric  and the difference tensor
   $K$ is defined as follows
 \begin{equation}
K(e_1,e_1)=e_2,\ \ \ K(e_1,e_2)=e_1,\ \ \  K(e_2,e_2)=3e_2,
 \end{equation}
where $e_1,e_2$ is the canonical basis of $\R^2$. The vector $U=e_1$
satisfies the equality in (\ref{nierownosc1a}). We also see that for
this structure the equality holds in (\ref{nierownosc2}).}
\end{example}

\begin{thm}
For any statistical structure, we have
\begin{equation}\label{nierownosc7}
\hat\rho-\rho=\Vert  A\Vert^2-\Vert E\Vert^2\ge -\frac{n-1}{3}\Vert
A\Vert^2,
\end{equation}
\begin{equation}\label{nierownosc8}
\hat\rho-\rho=\Vert  A\Vert^2-\Vert E\Vert^2\ge
-\frac{n-1}{n+2}\Vert E\Vert^2.
\end{equation}
For a statistical structure  on a Lagrangian submanifold of a
complex space form $\tilde M(4c)$, we have

\begin{equation}
\hat\rho\le \frac{n(n-1)}{n+2}\left(c(n+2)+n\Vert \mu\Vert^2\right),
\end{equation}
\begin{equation}
\hat\rho\le \frac{n-1}{3}\left(3cn+\Vert
A\Vert^2\right).\end{equation}
 For a statistical structure on an
equiaffine sphere with $R=HR_0$, we have

\begin{equation}
\hat\rho\ge \frac{n-1}{n+2}\left(Hn(n+2)-\Vert E\Vert^2\right),
\end{equation}

\begin{equation}
\hat\rho\ge \frac{n-1}{3}\left( 3Hn-\Vert A\Vert^2 \right).
\end{equation}
\end{thm}

\medskip
\section{Using Simons' formulas}

   We shall  use the following version of Simons' formula. Let
$g$ be a positive-definite metric tensor field on $M$. For any
tensor field $s$ on a manifold $M$ we have
\begin{equation}\label{Simons'}
\frac{1}{2}\Delta(\Vert s\Vert^2)=g(\Delta  s,s) +\Vert\hat\nabla
s\Vert^2,
\end{equation}
where the meaning of the Laplacians is explained in Preliminaries.
We shall consider tensor fields of type $(0,k)$. Recall the Ricci
identity:
\begin{equation}\label{Ricci}
(\hnabla^2s)(X,Y,...)-(\hnabla^2s)(Y,X,...)
 =(\hat R(X,Y)\cdot s) (...).
\end{equation}
In particular, if $\hat\nabla ^2s$ is symmetric, then $\hat
R(X,Y)\cdot s=0$ for every $X,Y$.
\medskip
We shall employ (\ref{Simons'}) for studying the statistical cubic
form $A$, but first, we shall make some comments on the applications
of (\ref{Simons'}) to $1$-forms and bilinear symmetric forms.

\subsection{For $1$-forms}
If $\tau$ is a $1$-form on $M$, then we have
\begin{equation}\label{Ricci1form}
(\hnabla^2\tau)(X,Y,Z)-(\hnabla^2\tau)(Y,X,Z)=-\tau(\hat
R(X,Y)Z)=g(\hat R(X,Y)E,Z),
\end{equation}
where $E=\tau^\sharp$. From the last formula we immediately get
\begin{proposition}
If $\tau$ is a closed $1$-form and $\hat R=0$, then the cubic form
$\hnabla^2\tau$ is symmetric.
\end{proposition}

\begin{proposition}\label{symmetrictau}
Assume that $\hnabla ^2\tau $ is a symmetric cubic form. Then $\im
\hat R\subset\ker\tau$. In particular, $ \widehat\Ric(X,E)=0 $ for
every $X$, hence $\tau=0$ if $\widehat\Ric$ is non-degenerate.
\end{proposition}

 By (\ref{Ricci1form}), we also get
\begin{equation}\label{Weitzenbockfor1form}
\tr_g(\hnabla^2 \tau)(\cdot,\cdot,X)=(d\delta \tau+\delta
d\tau)(X)+\widehat\Ric(X,E),
\end{equation}
Simons' formula (\ref{Simons'}) now yields

\begin{equation}\label{for1-form}
\frac{1}{2}\Delta (\Vert\tau\Vert ^2)=g((d\delta+\delta d)
\tau,\tau)+\widehat{\Ric}(E,E)+\Vert\hat\nabla\tau\Vert^2.
\end{equation}
 By Proposition \ref{symmetrictau}, (\ref{Weitzenbockfor1form}), and
 (\ref{for1-form}),
 we obtain
\begin{thm}\label{nablasquaretau}
Let $\tau$ be a $1$-form such that $\hnabla ^2\tau=0$. Then
$\widehat\Ric (X,E)=0$ for every $X$ and $\tau$ is harmonic.
\newline
If, additionally, $M$ is compact, or $\Vert \tau\Vert$ is constant,
then $\hat\nabla\tau=0$. If $\widehat\Ric$ is non-degenerate at a
point of $M$, then $\tau=0$.
\end{thm}

Of course, (\ref{for1-form})  immediately implies the following
classical Bochner's theorem
\begin{thm}\label{Bochner}
Let $(M, g)$ be  a compact  Riemannian manifold. Assume that
$\widehat{\Ric}$ is semi-positive definite on $M$. Then each
harmonic $1$-form $\tau$ on $M$ is $\hat\nabla$-parallel and
$\widehat\Ric (E,E)=0$ on $M$. If, moreover,  $\widehat\Ric$ is
positive definite at some point of $M$, then $\tau= 0$ on $M$.
\end{thm}

%2) $\dim M=2$, then either $\tau= 0$ or $\hat R=0$ on  $M$.
%3) $M$ is a topological sphere then $\tau=0$ on $M$. To z zaczesywania sfery

We can apply the above theorems to the Czebyshev form of a
statistical structure. For instance, we have

\begin{corollary}
Let $(g,\nabla)$ be a statistical structure on a manifold $M$ and
$\tau$ be its Czebyshev form. Assume that $\hnabla^2\tau =0$, $\Vert
\tau\Vert$ is constant or $M$ is compact and, moreover, the Ricci
tensor of the metric $g$ is non-degenerate at a point of $M$. Then
the structure is trace-free.
\end{corollary}

Recall that if $\dim M=2$, then  the Ricci tensor (at a point of
$M$) of a Levi-Civita connection is either non-degenerate or
vanishes. In the last case, the curvature tensor also vanishes at
the point. It provides additional information in the above theorems
in the case where $\dim M=2$. We have, for example,

\begin{corollary}
Let $\tau$ be a $1$-form on a $2$-dimensional manifold $M$ and
$\hnabla^2\tau=0$. At each point of $M$ we have $\tau=0$ or $\hat
R=0$. In the case where $M$ is compact, we have $\tau=0$ or $\hat
R=0$ on $M$.
\end{corollary}

%If $\div E=0$, where $E=JH$ and $H$ is the mean curvature vector of
%a compact Lagrangian submanifold of a complex space form, and the
%Ricci tensor of the submanifold is  semi-positive on $M$, then the
%mean curvature vector is parallel. If moreover the Ricci tensor is
%positive at least at one point of $M$, then the submanifold is
%minimal.
%\end{corollary}

\medskip
\subsection{For tensor fields of type $(0,2)$}

If $\beta$ is a symmetric bilinear form on a Riemannian manifold
$(M,g)$, then we have

$$
(\hat R(X,Y)\cdot \beta)(U,V)=-\beta(\hat R(X,Y)U,V)-\beta(U,\hat
R(X,Y)V).
$$
It follows that if $e_1,...,e_n$ is an orthonormal basis
diagonalizing $\beta$ and $\lambda_1,...,\lambda _n$ are the
corresponding eigenvalues, then
\begin{equation}
(\hat R(X,Y)\cdot \beta)(e_i,e_j)=g(\hat
R(X,Y)e_i,e_j)(\lambda_i-\lambda_j)=g(\hat
R(e_i,e_j)Y,X)(\lambda_i-\lambda_j).
\end{equation}
In particular, we have

\begin{proposition}
If $\hat\nabla^2\beta $ is symmetric and the sectional curvature of
$g$ is nonzero for every plane, then $\beta=\lambda g$.
\end{proposition}

 Simons' formula for symmetric tensor fields of type $(0,2)$ can be
formulated as follows
\begin{thm} Let $\beta$ be a symmetric tensor field of type $(0,2)$
on a Riemannian manifold $M$. If $\hat\nabla\beta$ is symmetric,
then
 we have
 \begin{equation}\label{for_type_(0,2)}
\begin{array}{rcl}
&&\frac{1}{2}\Delta\Vert\beta\Vert^2  =\Vert \hnabla\beta\Vert^2+
\sum_{ijk}\hat\nabla^2\beta(e_j,e_k, e_i,e_i)\beta_{jk}\\&&\ \ \ \ \
 \ \ \ \ \ \ \ +\sum_{i<k}\hat k(e_i\wedge e_k)(\lambda_i-\lambda_k)^2,
\end{array}
 \end{equation}
where $\beta_{jk}=\beta (e_j,e_k)$,
$\beta_{jk}=\delta_{jk}\lambda_k$ for some orthonormal basis
$e_1,...,e_n$ of $T_xM$, $x\in M$ and $\hat k(e_i\wedge e_k)$ is the
sectional curvature of $g$ by the plane spanned by $e_i, e_k$.
\end{thm}
\proof First we have
\begin{eqnarray*}
\sum_i\hnabla^2
\beta(e_i,e_i,X,Y)=\sum_i[\hnabla^2\beta(X,Y,e_i,e_i)+\beta(\hat
R(X,e_i))Y,e_i)+\beta(Y,\hat R(X,e_i)e_i)]
\end{eqnarray*}
for any orhonormal basis $e_1,...,e_n$ and any $X,Y$. Hence
\begin{eqnarray*}
&&g(\Delta\beta,\beta)=\sum_{kli}\hnabla^2\beta(e_k,e_l,
e_i,e_i)\beta_{kl}\\
&&\ \ \
 \ \ \ \ \ \ \ +\sum_{kils}\hat
 R_{kils}\beta_{is}\beta_{kl}+\sum_{ksl}\widehat{\Ric}_{ks}\beta_{kl}\beta_{ls}.
\end{eqnarray*}
If $e_1,..., e_n$ is a basis diagonalizing $\beta$ with eigenvalues
$\lambda_1,...,\lambda_n$, then we get
\begin{eqnarray*}
&&\sum_{kils} \hat
R_{kils}\beta_{is}\beta_{kl}+\sum_{ksl}\widehat{\Ric}_{ks}\beta_{kl}\beta_{ls}
 =-\sum_{ik}\hat R_{ikki}\lambda_i\lambda_k+\sum_k
 \widehat{\Ric}_{kk}\lambda_{k}^2\\
%&&=-2\sum_{i<k}\hat k(e_i\wedge e_k)\lambda_i\lambda_k+\sum_{i\ne
%k}\hat k(e_i\wedge e_k)\lambda _k^2\\
&&=\sum_{i<k}[-2\hat k(e_i\wedge e_k)\lambda_i\lambda_k+\hat
k(e_i\wedge e_k)\lambda_i^2 +\hat k(e_i\wedge e_k)\lambda_k^2]\\
&&\ \ \ \ \ \ \ \ \ = \sum_{i<k}\hat k(e_i\wedge
e_k)(\lambda_i-\lambda_k)^2.
\end{eqnarray*}
It is now sufficient to use (\ref{Simons'}).

\koniec

From the above formula, we immediately get

\begin{thm}\label{tw.dla2-form}
Let $(M,g)$ be a compact Riemannian manifold and $\beta$ be a
symmetric $(0,2)$-tensor field such that $\hat\nabla\beta$ is
symmetric and the function $\lambda =\tr_g\beta$ is constant on $M$.
If the sectional curvature of $g$ is  nonnegative on $M$ and
positive at some point of $M$, then $\beta=\frac{\lambda}{n}g$.
\end{thm}

\proof If the sectional curvature of $g$ is everywhere nonnegative,
 then by  (\ref{for_type_(0,2)}) we know that $\hnabla \beta=0$.
 The formula (\ref{for_type_(0,2)}) implies also that for every
 index $i\ne k$ $\hat k(e_i\wedge e_k)(\lambda _i-\lambda
 _k)^2=0$. If at some point of $M$ all sectional curvatures $\hat k$
 are positive, we see that $ \beta =\frac{\lambda}{n}g$ at this point and
 hence everywhere on $M$.\koniec

In the above theorem, instead of assuming that $M$ is compact, one
can assume that $\Vert\beta\Vert$ is constant on $M$.

\begin{thm} let $(M,g)$ be a  connected Riemannian  manifold  whose sectional curvature
is   positive (or negative) at some point of $M$. If $\beta$ is a
symmetric 2-form and $\hnabla\beta=0$, then $\beta=c g$ for some
constant $c$.
\end{thm}

\begin{corollary} Let $(M,g)$ be a compact Riemannian manifold whose
sectional curvature is nonnegative on $M$ and positive at some point
of $M$. If $\hnabla \widehat\Ric $ is symmetric and the scalar
curvature $\hat\rho $ is constant, then the Riemannian structure is
Einstein.
\end{corollary}

\begin{corollary} Let $(M,g)$ be a Riemannian manifold whose
sectional curvature is  positive at some point of $M$ (or negative
at some point of $M$). If $\hnabla \widehat\Ric =0$, then the
Riemannian structure is Einstein.
\end{corollary}

%\begin{corollary}
%Let $(M,g,\nabla)$ be a Hessian manifold whose metric tensor field
%$g$ has non-negative sectional curvature on $M$ and positive
%sectional curvature at some point of $M$. Let $\beta=\nabla \tau$ be
%the Koszul form. If $\hnabla \beta$ is symmetric, $\Vert\beta\Vert$
%is constant, $\tr_g\beta$ is constant (i.e. $d\delta\tau
%=d(\Vert\tau\Vert^2)$), then the structure $(g,\nabla)$ is
%Einstein-Hessian.
%\end{corollary}
In the geometry of statistical structures, we have a few symmetric
bilinear forms, for instance,  $\tau\circ K$, $g(K_\cdot,K_\cdot)$.
The Czebyshev form is closed  if and only if the bilinear  forms
$\nabla \tau$, $\Ric$, $\overline\Ric$ are also symmetric. In
particular, we can formulate the following

\begin{corollary}
Let $(M,g,\nabla)$ be a connected statistical manifold. If $\Ric$ is
symmetric, $\hnabla \Ric=0$ and the sectional curvature of $g$ is
positive (or negative) at some point of $M$, then $\Ric=\lambda g$
for some constant $\lambda$.
\end{corollary}

\begin{corollary} Let $(M,g,\nabla)$ be a Hessian manifold and
$\beta=\nabla \tau$ its Koszul form. If $\hnabla \beta=0$ and the
sectional curvature for $g$ is positive (or negative) at some point
of $M$, then the structure is Einstein-Hessian.
\end{corollary}

\medskip

\subsection{For cubic forms}
We shall first compute a Simons' type formula for symmetric cubic
forms on an $n$-dimensional  Riemannian manifold $(M,g)$. Let $A$ be
a symmetric cubic form and $K=A^\sharp$. Assume that $\hnabla A$ is
symmetric. Thus we have a conjugate symmetric statistical structure.
Let $e_1,...,e_n$ be an orthonormal basis of a tangent space $T_xM$.
We have
\begin{eqnarray*}
\begin{array}{rcl}
&&\hnabla^2A(e_i,e_i, X,Y,Z)=\hnabla^2A(X,Y,Z, e_i,e_i)-A(\hat
R(e_i,X)Y, Z,e_i)\\
&&\ \ \ \ \ \ \ \  \ \ \ \ \ \ \ \ \ -A(Y,\hat
R(e_i,X)Z,e_i)-A(Y,Z,\hat R(e_i,X)e_i)
\end{array}
\end{eqnarray*}
for every $i$. Hence

\begin{eqnarray*}
\begin{array}{rcl}
&&\hnabla A(e_i,e_i, e_j,e_k, e_l)=\hnabla
A(e_j,e_k,e_l,e_i,e_i)-g(K(e_l,e_i), \hat R(e_i,e_j)e_k)\\
&&\ \ \ \ \ \ \ \ \ \ \ \ \ \ \ \ \ \ \ \ \  -g(K(e_k,e_i),\hat
R(e_i,e_j)e_l)-g(K(e_k,e_l),\hat R(e_i,e_j)e_i)\\
&&\ \ \ \ \ \ \ =\hnabla A(e_j,e_k,e_l,e_i,e_i)-\sum _r A_{lir}\hat
R_{ijkr}-\sum _r A_{kir}\hat R_{ijlr}-\sum_rA_{klr}\hat R_{ijir}
\end{array}
\end{eqnarray*}
for every $i,j,k,l$.
 We also have
\begin{equation}
\hnabla ^2\tau(X,Y,Z)=\sum _i\hnabla^2 A(X,Y,Z, e_i,e_i).
\end{equation}
We now compute

\begin{eqnarray*}
\begin{array}{rcl}
&&\sum_{ijkl}\hnabla ^2A(e_i,e_i, e_j, e_k,
e_l)A_{jkl}=\sum_{ijkl}\hnabla ^2A(e_j,e_k, e_l, e_i, e_i)A_{jkl}\\
&&\ \ \ \ \ \ -\sum_{ijklr}[A_{lir}A_{jkl}\hat R_{ijkr}+
A_{kir}A_{jkl}\hat R_{ijlr}+A_{klr}A_{jkl}\hat R_{ijir}]\\
&&\ \ \ \ \ \ =g(\hnabla^2 \tau,
A)+\sum_{kljr}\widehat{\Ric}_{jr}A_{jkl}A_{rkl}\\
&&\ \ \ \ \ \ \ \ \ \ \ \ \ \ \ \ \ \ \ \ \ \ \ \ \ \ \ \ \ \ \ \ \
\ \ \ -\sum_{ijklr} A_{irl}A_{jkl}\hat
R_{ijkr}-\sum_{ijklr}A_{irk}A_{jlk}\hat
R_{ijlr}\\
&&\ \ \ \ \ \ =g(\hnabla^2 \tau,
A)+\sum_{kljr}\widehat{\Ric}_{jr}A_{jkl}A_{rkl}+\sum_{ijklr}(A_{ikr}A_{jlr}-A_{ilr}A_{jkr})\hat
R_{ijkl}\\
&&\ \ \ \ \ \ =g(\hnabla^2 \tau, A)+ g(\widehat {\Ric}, g(K_\cdot ,
K_\cdot))-\sum_{ijkl}g([K_{e_i}, K_{e_j}]e_k,e_l)\hat R_{ijkl}.
\end{array}
\end{eqnarray*}

Since for a conjugate symmetric statistical structure $R=\hat
R+[K,K]$ and $\Ric=\widehat\Ric+\tau\circ K-g(K_\cdot,K_\cdot )$
(see Preliminaries), by using (\ref{Simons'}), we obtain

\begin{thm}
For any conjugate symmetric statistical structure, we have
\begin{equation}\label{1dla A}
\frac{1}{2}\Delta (\Vert A\Vert^2)=\Vert\hnabla A\Vert^2+
g(\hnabla^2 \tau ,A)-g([K,K],\hat
R)+g(\widehat{\Ric},g(K_\cdot,K_\cdot)),
\end{equation}
\begin{equation}\label{2dla A}
\frac{1}{2}\Delta (\Vert A\Vert^2)=\Vert\hnabla A\Vert^2+
g(\hnabla^2 \tau ,A)+g(\hat R-R,\hat
R)+g(\widehat{\Ric},g(K_\cdot,K_\cdot)),
\end{equation}

\begin{equation}\label{3dla A}
\begin{array}{rcl}
&&\frac{1}{2}\Delta (\Vert A\Vert^2)=\Vert\hnabla A\Vert^2+
g(\hnabla^2 \tau ,A)\\
&&\ \ \ \ \ \ \ \  \ \ \ \ \ \ \ \ \ \ \ \  +\hat
R^2+\widehat{\Ric}^2-g(R,\hat
R)-g(\Ric,\widehat{\Ric})\\
&&\ \ \ \ \ \ \ \ \ \ \ \ \ \ \ \ \ \ \  \ \ \  \ \ \ \ \
 \ \ \  \ \ \ \ \ \ \ \ \ \ \ \ \ \ \ \ \ \ \ \ \ \ \ \ \ \ \ \ \ \
 +g(\widehat{\Ric}, \tau\circ K).
\end{array}
\end{equation}
\end{thm}

Consider now the case where the sectional $K$-curvature is constant.
In \cite{seccur} the following result was proved

\begin{thm}
Let  $(g,K)$ be a conjugate symmetric trace-free statistical
structure on a manifold $M$. If the sectional $K$-curvature is
constant (automatically nonpositive, because of the trace-freeness),
then either $K=0$ or: $\hat R=0$ and $\hat\nabla K=0$.
\end{thm}

In a more general situation, we get
\begin{thm} Let $(g, K)$ be a conjugate symmetric statistical
structure with $\hat\nabla$-parallel Czebyshev form $\tau$ and
$\widehat{\Ric}\ge 0$. If the sectional $K$-curvature is constant
and nonpositive, then $\hat\nabla K=0$. If the sectional
$K$-curvature is a negative constant,  then  $\widehat{\Ric}=0$.
\end{thm}
\proof If $[K,K]=\kappa R_0$, then (\ref{1dla A}) becomes

\begin{equation}\label{4dla A}
\frac{1}{2}\Delta (\Vert A\Vert^2)=\Vert\hnabla A\Vert^2+
g(\hnabla^2 \tau ,A)-2\kappa\hat\rho
+g(\widehat{\Ric},g(K_\cdot,K_\cdot)).
\end{equation}
Since $\rho^K=\Vert E\Vert^2-\Vert A\Vert^2$ (see (\ref{rho^K})) and
$\hat\nabla E=0$, we have that $\Vert A\Vert$ is constant. Since
$\widehat{\Ric}\ge 0$, we have
$g(\widehat{\Ric},g(K_\cdot,K_\cdot))\ge 0$. The first assertion now
follows from (\ref{4dla A}). If $\kappa< 0$, then $\hat\rho=0$.
Since $\widehat{\Ric}\ge 0$, we have $\widehat{\Ric}=0$.
 \koniec

For trace-free statistical structures (\ref{3dla A}) becomes

\begin{equation}\label{3a_dla A}
\frac{1}{2}\Delta (\Vert A\Vert^2)=\Vert\hnabla A\Vert^2+\hat
R^2+\widehat{\Ric}^2-g(R,\hat R)-g(\Ric,\widehat{\Ric}).
\end{equation}

 The last formula was proved by
An-Min Li for hyperbolic affine Blaschke spheres, \cite{L}.
 In that case, $\tau =0$, $R=HR_0$, where $H$ is a constant
  and $\Ric=(n-1)Hg$. Hence (\ref{3a_dla A}) can be written as
  follows

\begin{equation}\label{5dla A}
\frac{1}{2}\Delta (\Vert A\Vert^2)=\Vert\hnabla A\Vert^2+ \hat
R^2+\widehat\Ric^2-(n+1)H\hat\rho.
\end{equation}
Moreover,  $\Vert A\Vert^2= \hat\rho -n(n-1)H$, see
(\ref{theorema-egregium}). Hence $\Vert A\Vert$ is constant if
$\hat\rho$ is constant. If $\hat\rho=0$, $H$ must be nonpositive.
Thus we have
\begin{thm}{\cite{LSZ}}
For a  Blaschke affine sphere if $\hat\rho=0$, then $\hat R =0$ and
$\hnabla A=0$. The sphere is either a part of a trivial paraboloid
(an elliptic paraboloid with its trivial affine structure) or the
sphere is hyperbolic.
\end{thm}
The hyperbolic sphere from the last theorem must be   given by the
equation $x_1\cdot ...\cdot x_n=c$, where $c$ is a positive
constant, see \cite{LSZ}.

For equiaffine spheres, we  can prove

\begin{thm} For an equiaffine hyperbolic or parabolic sphere, if $\hnabla
\tau=0$,  $\widehat{\Ric}\ge0$  and $\hat\rho $ is constant, then
$\hat R=0$ and $\hnabla A=0$.\end{thm}

\proof Since $R=HR_0$ and $\hat\nabla \tau=0$, the formula
(\ref{2dla A}) can be written as follows
\begin{equation}\label{6dla A}
\frac{1}{2}\Delta (\Vert A\Vert^2)=\Vert\hnabla A\Vert^2+\hat
R^2-2H\hat\rho +g(\widehat{\Ric},g(K_\cdot,K_\cdot)).
\end{equation}
Moreover, the function $\Vert E\Vert$ is constant. Since the
functions $\hat\rho$ and  $\rho $ are constant,
(\ref{theorema-egregium}) implies that $\Vert A\Vert$ is also
constant. The assertion now follows from (\ref{6dla A}).\koniec

%Similarly one gets To wynika z czesci o 1-formach Let $(g,\nabla)$
%be a conjugate symmetric statistical structure on a compact manifold
%$M$. If the sectional $\nabla$-curvature  is constant non-positive,
%$\widehat\Ric\ge 0$ and $\hnabla^2\tau=0$, then $\hat R=0$ and
%$\hnabla A=0$.
%\end{thm}

From (\ref{1dla A}), we immediately get

\begin{proposition}
Let $(M,g,\nabla)$ be a conjugate symmetric statistical manifold
such that $\hnabla ^2\tau=0$, $\hat R=0$, and $\Vert A\Vert$ is
constant (or $M$ is compact). Then $\hnabla A=0$.
\end{proposition}
From (\ref{6dla A}) one gets

\begin{proposition}
For a Hessian structure $(g,A)$, if $\hnabla A=0$  and
$\widehat\Ric\ge 0$, then $\hat R=0$.

\end{proposition}

\bigskip
For a Lagrangian submanifold in a complex space form $\tilde M(4c)$,
we have $ -[K,K]=cR_0-\hat R $, see (\ref{gausslagrangian}). Using
now (\ref{1dla A})  we obtain
%\begin{equation}
%\Delta\Vert A\Vert^2=\Vert\hat\nabla A\Vert^2-\hat R^2-\widehat
%\Ric^2+c(n+1)\hat\rho +g(\widehat\Ric,\tau\circ K ).
%\end{equation}

\begin{equation}\label{forLagrangian}
\Delta\Vert A\Vert^2=\Vert\hat\nabla A\Vert^2
+g(\hat\nabla^2\tau,A)-\Vert\hat R\Vert^2+2c\hat\rho
+g(\widehat\Ric,g(K_\cdot, K_\cdot)).
\end{equation}

\begin{thm} Let $(g,K)$ be the statistical structure on a Lagrangian
submanifold in  a complex space form $\tilde M (4c)$ whose second
fundamental tensor $K$ is $\hat\nabla$-parallel. If
$\widehat\Ric=0$, then $\hat R=0$. If  $\widehat\Ric\le 0$ and $c\ge
0$ then $\hat R=0$.
\end{thm}

If $g$ has constant curvature $\hat k$, then using
(\ref{theorema_egregium_Lagrangian}) and (\ref{forLagrangian}) we
get

\begin{equation}
\Delta\Vert A\Vert^2=\Vert\hat\nabla A\Vert^2
+g(\hat\nabla^2\tau,A)+\hat k[2(\Vert A\Vert ^2-\Vert E\Vert
^2)+(n-1)\Vert A\Vert^2].
\end{equation}

Hence, if $\hat k\ge 0$, then (\ref{nierownosc8}) yields
\begin{equation}\label{forLagrangian1}
\Delta\Vert A\Vert^2\ge\Vert\hat\nabla A\Vert^2
+g(\hat\nabla^2\tau,A)+\frac{n-1}{3}\hat k\Vert A\Vert^2.
\end{equation}

This formula immediately yields a theorem proved by Chen and Ogiue
in \cite{ChBY46} and saying that a minimal Lagrangian submanifold of
constant positive curvature in a complex space form must be totally
geodesic. Indeed, it is sufficient to note that if $\tau=0$, then,
by (\ref{theorema_egregium_Lagrangian}), $\Vert A\Vert$ is constant.
Using (\ref{theorema_egregium_Lagrangian}) and Theorem
\ref{nablasquaretau}, we also get
\begin{thm}
Let $M$ be a  flat Lagrangian submanifold in a complex space form
and  $\hat\nabla ^2\tau=0$. If $M$ is compact or the mean curvature
 is constant, then the submanifold has the parallel second fundamental
tensor.
\end{thm}

 \bigskip

\section{Using a maximum principle} In this section, we set $u=\Vert A\Vert ^2$ and we assume  that
$(g,\nabla)$ is a trace-free statistical structure on an
$n$-dimensional connected manifold $M$ such that $R=HR_0$. It is
automatically  conjugate symmetric and therefore, if $n>2$, then $H$
is constant. If $n=2$,  $H$  can be a function. In this case, the
assumption $R=HR_0$ is equivalent to the assumption that
$(g,\nabla)$ is conjugate symmetric. If $H$ is constant, such a
statistical structure can be automatically realized (locally) on an
affine Blaschke sphere.

We shall now employ two inequalities due to Calabi and An-Min Li.
For detailed proofs, we refer to  \cite{NS} and \cite{L_1}. Here we
give a sketch of the  proof. First define  tensor fields $L$, $P$
and $Q$ as follows
\begin{equation}
L(X,Y,W,Z)=g(K(X,Y),K(W,Z)),
\end{equation}

\begin{equation}
P(X,Y,W,Z)=L(X,Y,W,Z)-L(W,Y,X,Z),
 \end{equation}
\begin{equation}
Q(Y,W,Z)=\tr_g([K_{\cdot}, K_Y]\cdot A)(\cdot,W,Z).
\end{equation}
In the case under consideration,  we have $\hat R=HR_0-[K,K]$ and
$$\Delta A (Y,W,Z)=\tr_g(\hat R(\cdot,Y)\cdot A)(\cdot,W,Z)=H(n+1)A
-Q(Y,W,Z).$$ We have
\begin{equation}
\Vert L\Vert ^2=\sum_{ij}a_{ij}^2,\ \ \ \ \ \ \  \Vert
P\Vert^2=\sum_{ijkl}b_{ij:kl}^2,
\end{equation}
where
\begin{equation}\label{notationNS}
\begin{array}{rcl}
&&a_{ij}=\sum_{kl}A_{ikl}A_{jkl},\\
&&b_{ij;kl}=\sum_m(A_{ijm}A_{klm}-A_{kjm}A_{ilm}),
\end{array}
\end{equation}
The notation in (\ref{notationNS}) is similar to that in \cite{NS}
on p. 84, where it was introduced for the cubic form $C=-2A$.
 It was proved
there that $-g(Q,A)=\Vert L\Vert^2+\Vert P\Vert^2$, that is,

\begin{equation}
g(\Delta A,A)=(n+1)Hu+\Vert L\Vert^2+\Vert P\Vert^2.
\end{equation}
Moreover, see \cite{NS},
\begin{equation}\label{L^2+P^2_1}
\Vert L\Vert^2+\Vert P\Vert^2\ge \frac{n+1}{n(n-1)}u^2,
\end{equation}
where $u=\Vert A\Vert^2$. In  \cite{NS} it was assumed that $H$ is
constant, but this assumption is  not necessary (although
automatically satisfied if $n>2$). The inequality  actually is due
to Blaschke (for $n=2$) and Calabi (for $n>2$). In \cite{L_1} A-M.
Li proved the following algebraic inequality (Theorem 1   in
\cite{L_1})

\begin{equation}\label{L^2+P^2_2}
\Vert L\Vert^2+\Vert P\Vert^2\le \frac{3}{2}u^2.
\end{equation}
In particular, in the case where $n=2$, by comparing
(\ref{L^2+P^2_1}) and (\ref{L^2+P^2_2}), we get
\begin{equation}
\Vert L\Vert^2+\Vert P\Vert^2=\frac{3}{2}u^2.
\end{equation}
The last equality was also directly proved in \cite{NS}. Using the
above facts and Simons' formula (\ref{Simons'}), one can formulate
\begin{thm}\label{tw.NSL}
Let $(g,\nabla)$ be a trace-free statistical structure on an
n-dimensional manifold $M$ and $R=HR_0$. We have
\begin{equation}\label{inequalitiesNSLI}
(n+1)Hu +\frac{n+1}{n(n-1)}u^2+\Vert\hat\nabla
A\Vert^2\le\frac{1}{2}\Delta u\le
(n+1)Hu+\frac{3}{2}u^2+\Vert\hat\nabla A\Vert^2.
\end{equation}
If $n=2$, the equalities hold.
\end{thm}

As an immediate consequence of this theorem, we have
\begin{corollary}
Let $(g,\nabla)$ be a trace-free statistical structure on an
n-dimensional manifold $M$ and $R=HR_0$. Assume that $\hat\nabla
A=0$ on $M$.
\newline
1) If $H\ge 0$, then $A=0$ on $M$.
\newline
2) If $H< 0$, then either $A=0$ on $M$ or
\begin{equation}
\frac{2}{3}(n+1)(-H)\le u\le n(n-1)(-H).
\end{equation}
\newline
If $n=2$ and $\sup\, H=0$, then $A=0$.
\end{corollary}
For  treating the case where $g$ is complete, we shall use the
following maximum principle

\begin{lemma}\label{Y}\cite{Y}.
Let $(M,g)$ be a complete Riemannian manifold, whose Ricci tensor is
bounded from below and let $f$ be a $\mathcal C^2$-function bounded
from above on $M$.  For every $\varepsilon >0$, there exists a point
$x\in M$  at which
\newline
i) $f(x)> \sup f -\varepsilon$
\newline
ii)$\Delta f(x)<\varepsilon$.

\end{lemma}

The results contained in the following theorem and their
generalizations can be found in \cite{Cal}, \cite{LSZ},
\cite{calabi}.

\begin{thm}\label{realcalabi} Let $(g,\nabla)$ be a trace-free
statistical structure with $R=HR_0$ and complete $g$ on a manifold
$M$. If $H\ge 0$, then the structure is trivial. If $H$ is constant
and negative, then the Ricci tensor of $g$ is nonpositive and
consequently
\begin{equation}\label{Calestimation}
u\le n(n-1)(-H).
\end{equation}
\end{thm}
We shall now give a more delicate estimation for $u$ by looking more
carefully at  the function $\Vert \hnabla A\Vert$. Recall here that
for a trace-free conjugate symmetric statistical structure we have
$\widehat\Ric>\Ric$, see (\ref{Riccihat_ge}). Hence, in all theorems
below, the Ricci tensor of $g$ is automatically bounded from below.

\begin{thm}\label{aboutinfu}
Let $(g,\nabla)$ be a trace-free statistical structure on an
$n$-dimensional  manifold $M$ and $R=HR_0$, where $H$ is a negative
constant. If $g$ is complete  and
\begin{equation}
\sup\, \Vert\hat\nabla A\Vert^2<\frac{H^2(n+1)^2}{6},
\end{equation} then
\begin{equation}
\inf u\ge \frac{(n+1)(-H)+\sqrt{(n+1)^2H^2-6N_2}}{3},
\end{equation}
or
\begin{equation}
\inf u\le \frac{(n+1)(-H)-\sqrt{(n+1)^2H^2-6N_2}}{3},
\end{equation}
 where  $u=\Vert A\Vert^2$ and $N_2=\sup \Vert \hat\nabla
A\Vert^2$.
\end{thm}
\proof By Theorem \ref{tw.NSL}, we have
\begin{equation}\label{1dlanablaA}
\frac{1}{2}\Delta u\le (n+1)Hu+\frac{3}{2}u^2+N_2.
\end{equation}
Set $N_1=\inf u$ and $\tilde u=-u$. The function $\tilde u$ is
bounded from above by $0$ and $-N_1=\sup\, \tilde u$. Take sequences
$\varepsilon _i$ and $x_i$ satisfying $i)$ from  Lemma \ref{Y}
applied to the function $\tilde u$. Assume that the sequence
$\varepsilon _i$ is decreasing and tending to $0$. We have $\tilde u
(x_i)>-N_1-\varepsilon_i>-N_1-\varepsilon _1$ and consequently
$u(x_i)<N_1+\varepsilon _1$ and $-u^2(x_i)>-(N_1+\varepsilon_1)^2$.
By (\ref{1dlanablaA}),  we now have
\begin{equation}\label{deltatildeu}
\frac{1}{2}(\Delta \tilde u)(x_i)\ge
-(n+1)HN_1-\frac{3}{2}(N_1+\varepsilon_1)^2-N_2.
\end{equation}
Consider the following polynomial of degree 2
\begin{equation}\label{polynomial}
\begin{array}{rcl}
&&p(t)=-\frac{3}{2}(t+\varepsilon _1)^2-(n+1)Ht-N_2\\
&&=
-\frac{3}{2}t^2+[-3\varepsilon_1-(n+1)H]t-[\frac{3}{2}\varepsilon_1^2+N_2].
\end{array}
\end{equation}
Denote by $\delta$ its  discriminant. We have
\begin{equation}\label{discriminat}
\delta=N^2+6\varepsilon_1(n+1)H<N^2,\end{equation} where
$N=\sqrt{(n+1)^2H^2-6N_2}$. Assume that $\varepsilon_1$ is so small
that $\delta>0$. Then $\delta=N^2q^2$ for some $0<q<1$. The relation
between $q$ and $\varepsilon_1$ is the following

\begin{equation}
\varepsilon_1=\frac{N^2(1-q^2)}{6(n+1)(-H)}.
\end{equation}
Instead of choosing $\varepsilon_1$, we can choose $q$, and
$\varepsilon_1$ tends to $0$ if and only if $q$ tends to $1$. Let
$t_1$, $t_2$ be the roots of the polynomial $p(t)$ depending also on
$q$. We have

\begin{equation}
\begin{array}{rcl}
&&t_1(q)=\frac{(n+1)(-H)-q\sqrt{(n+1)^2H^2-6N_2}}{3}-\frac{N^2(1-q^2)}{6(n+1)(-H)}\\
&& t_2(q)
=\frac{(n+1)(-H)+q\sqrt{(n+1)^2H^2-6N_2}}{3}-\frac{N^2(1-q^2)}{6(n+1)(-H)}
\end{array}
\end{equation}
One sees that $$t_1(q)\to
\frac{(n+1)(-H)-\sqrt{(n+1)^2H^2-6N_2}}{3},$$
$$t_2(q)\to\frac{(n+1)(-H)+\sqrt{(n+1)^2H^2-6N_2}}{3}$$ if $q$
tends to $1$. We shall now use the maximum principle. Suppose  that
\begin{equation}
\begin{array}{rcl}
&&\frac{(n+1)(-H)-q\sqrt{(n+1)^2H^2-6N_2}}{3}<N_1\\
&&\ \ \ \ \ \ \ \ \ \ \ \ \ \ \ \ \ \ \ \ \ \ \ \
<\frac{(n+1)(-H)+q\sqrt{(n+1)^2H^2-6N_2}}{3}.
\end{array}
\end{equation} There is $q$ (sufficiently close to $1$) such that $t_1(q)<N_1<t_2(q)$. Hence
$p(N_1)$, depending only on $n$, $H$, $N_1$, $N_2$, and $\varepsilon
_1$ (equivalently $q$),  is positive and, by (\ref{deltatildeu}), we
cannot have $\frac{1}{2}\Delta\tilde u(x_i)< \varepsilon_i$ for
$\varepsilon_i$ tending to $0$. We have got a contradiction.\koniec

\begin{corollary}
Let $(g,\nabla)$ be a trace-free statistical structure on an
$n$-dimensional  manifold $M$ and $R=HR_0$, where $H$ is a negative
number. If $g$ is complete,
\begin{equation}
N_2<\frac{H^2(n+1)^2}{6} \end{equation} and
\begin{equation}\inf u> \frac{(n+1)(-H)-\sqrt{(n+1)^2H^2-6N_2}}{3},
\end{equation}
then
\begin{equation}
u\ge\frac{(n+1)(-H)+\sqrt{(n+1)^2H^2-6N_2}}{3}
\end{equation}
on $M$.
\end{corollary}

\begin{thm}\label{ostatnie}
Let $(g,\nabla)$ be a trace-free statistical structure on an
$n$-dimensional  manifold $M$ and $R=HR_0$, where $ H$ is a negative
number. If $g$ is complete, then
\begin{equation}\label{infN_4}
\inf\, \Vert\hat\nabla A\Vert^2\le\frac{n(n^2-1)H^2}{4}
\end{equation} and the following inequalities hold
\begin{equation}
\sup\, u\le \frac{n(n-1)(-H)+\sqrt{n^2(n-1)^2H^2-4N_4}}{2},
\end{equation}
\begin{equation}
\sup\, u\ge \frac{n(n-1)(-H)-\sqrt{n^2(n-1)^2H^2-4N_4}}{2},
\end{equation}
 where  $u=\Vert A\Vert^2$ and $N_4=\frac{n(n-1)}{n+1}\inf\, \Vert \hat\nabla
A\Vert^2$.
\end{thm}
\proof By Theorem \ref{tw.NSL}, we have
\begin{equation}\label{2dlanablaA}
\frac{1}{2}\Delta u\ge (n+1)H u
+\frac{n+1}{n(n-1)}u^2+\Vert\hat\nabla A\Vert^2.\end{equation} Let
$N_3=\sup\, u$. By Theorem \ref{realcalabi},  $N_3$ is finite. Take
sequences $\varepsilon _i$ and $x_i$ satisfying $i)$ from Lemma
\ref{Y} applied to the function $u$. Assume that $\varepsilon _i$ is
decreasing and tending to $0$, and $N_3-\varepsilon_1>0$. By
(\ref{2dlanablaA}), we have

\begin{equation}\label{deltau1}
\frac{n(n-1)}{2(n+1)}\Delta u(x_i)\ge n(n-1)HN_3 +
(N_3-\varepsilon_1)^2+N_4.\end{equation}
 Consider the following
polynomial of degree 2
\begin{equation}\label{polynomial1}
\begin{array}{rcl}
&&p(t)=-n(n-1)(-H)t+(t-\varepsilon_1)^2+N_4\\
&&\ \ \ \ \ \ \
=t^2-[2\varepsilon_1+n(n-1)(-H)]t+(\varepsilon_1^2+N_4).
\end{array}
\end{equation}
Denote by $\delta$ its  discriminant. We have
\begin{equation}\label{discriminat1}
\delta=n^2(n-1)^2H^2-4N_4+4n(n-1)(-H)\varepsilon_1.\end{equation} If
$n^2(n-1)^2H^2-4N_4<0$, then, by taking $\varepsilon_1$ sufficiently
small, we get $\delta<0$ and hence $\Delta u(x_i)$ will be positive
and bounded from zero by a positive number. This will give  a
contradiction to the maximum principle. We have proved
(\ref{infN_4}), that is, $n^2(n-1)^2H^2-4N_4\ge 0$.

Set $N^2=n^2(n-1)^2H^2-4N_4$. Assume first that  $N=0$.  Then
$\delta$ is positive and the roots of the polynomial
(\ref{polynomial1})  are
\begin{equation}
\frac{n(n-1)(-H)}{2}+\varepsilon_1\pm\sqrt{n(n-1)(-H)\varepsilon_1}.
\end{equation}
If $\varepsilon_1$ tends to $0$, both roots tend to
$\frac{n(n-1)(-H)}{2}$. Suppose that $N_3>\frac{n(n-1)(-H)}{2}$.
There is $\varepsilon_1$ such that
$N_3>\frac{n(n-1)(-H)}{2}+\varepsilon_1\pm\sqrt{n(n-1)(-H)\varepsilon_1}$
and then all $\Delta u(x_i)$ are positive and bounded from zero by a
positive number. This gives a contradiction with the maximum
principle. Similar arguments work in the case where
$N_3<\frac{n(n-1)(-H)}{2}$.

Assume now that $N^2>0$. There is $Q>1$ such that $\delta=N^2Q^2$.
The relation between $\varepsilon_1$ and $Q$ is the following

\begin{equation}
\varepsilon_1=\frac{N^2(Q^2-1)}{4n(n-1)(-H)}.
\end{equation}
The roots of  the polynomial (\ref{polynomial1})  are
\begin{equation}\label{rootspolynomial1}
\begin{array}{rcl}
&&t_1(Q)=\frac{n(n-1)(-H)(Q^2+1)-2Q\sqrt{n^2(n-1)^2H^2-4N_4}}{4}-\frac{N_4(Q^2-1)}{n(n-1)(-H)},\\
&&t_2(Q)=\frac{n(n-1)(-H)(Q^2+1)+2Q\sqrt{n^2(n-1)^2H^2-4N_4}}{4}-\frac{N_4(Q^2-1)}{n(n-1)(-H)}.
\end{array}
\end{equation}
If $Q\to 1$, then $$t_1(Q)\to
L_1=\frac{n(n-1)(-H)-\sqrt{n^2(n-1)^2H^2-4N_4}}{2}\ge 0$$ and
$$t_2(Q)\to P_1=\frac{n(n-1)(-H)+\sqrt{n^2(n-1)^2H^2-4N_4}}{2}\le
n(n-1)(-H).$$ Suppose now that  $N_3<L_1$ or $N_3>P_1$. If
$N_3<L_1$, then there is $Q$ (equivalently   the  corresponding
$\varepsilon_1)$ such that $N_3<t_1(Q)$, which means that $p(N_3)$
is positive and depends only on $n$, $H$, $N_4$ and $Q$. Hence
(\ref{deltau1}) gives a contradiction with the maximum principle. We
argue similarly if $N_3>P_1$. \koniec

The case where $n=2$ should be treated separately because in this
case $H$ may be a function. Note that $\Ric =Hg$, hence $\widehat
\Ric$ is bounded from below, if $H$ is bounded from below. The same
arguments as those used in the the above proofs yield
\begin{thm}\label{najostatnie}
 Let $(g,\nabla)$ be a  conjugate symmetric  trace-free statistical
structure with complete $g$ on a 2-dimensional manifold $M$. Then
$R=HR_0$ for some function $H$ and if $-\infty<H_2\le H\le 0$ for
some number $H_2$, then

\begin{equation}
\inf \hat u\le \frac{3}{2}H_2^2
\end{equation}
and
\begin{equation}
-H_2-\sqrt{H_2^2-\frac{2}{3}\inf\hat u}\le \sup \, u\le
-H_2+\sqrt{H_2^2-\frac{2}{3} \inf \hat u},
\end{equation}
where $\hat u=\Vert \hat \nabla A\Vert^2$.
 \newline
 If
$-\infty<H_2\le H\le H_1\le 0$ for some numbers $H_1$, $H_2$ and
\begin{equation}
\sup \,\hat u\le \frac{3}{2}H_1^2,
\end{equation}
then

\begin{equation}
\inf\, u\ge -H_1+\sqrt{H_1^2-\frac{2}{3}\sup \,\hat u}
\end{equation}
or
\begin{equation}
\inf\, u\le -H_1-\sqrt{H_1^2-\frac{2}{3}\sup\,\hat u}.
\end{equation}

\end{thm}

\bigskip

\section{Using Ros' integral formula}
 Let $s$ be a tensor of type $(0,k)$, $k\ge 2$ and  $g_0$ --
  the standard scalar product on $\R^n$.  Define the following $1$-form on $S^{n-1}=\{V\in \R^n;\Vert
V\Vert=1\}$:
$$\alpha_V(e)=s(V,...,V,e,V,...,V),$$ where $e$ stays at a fixed place $i_0$.
Let $\delta$ denote the codifferential relative to $g_0$ restricted
to $S^{n-1}$. By a straightforward computation, one gets
\begin{equation}\label{spherical_codifferential}
\begin{array}{rcl}
&&\delta \alpha =-(n+k-2) s(V,...,V)\\
&&\ \ \ \ \ \ \ +\tr_gs(\cdot,V,...,V,
\cdot,V,...,V)+...+\tr_gs(V,...,V,\cdot, V,...,V, \cdot),
\end{array}
\end{equation}
where one of the dots "$\cdot$ " stays at the fixed $i_0$-th place.

Let now  $g$ be a positive-definite metric tensor field on a
manifold $M$. For simplicity, we shall assume that $M$ is connected
and oriented. Let  $UM$ denote the unit sphere bundle over $M$. By
(\ref{spherical_codifferential}), we have

\begin{proposition} For a covariant tensor field $s$ of degree $k\ge 2$ we
have
\begin{equation}\label{integral_over_fiber}
\begin{array}{rcl}
&&(n+k-2)\int_{U_xM} s(V,...,V)\\
&&\ \ \ \ \ \ \ \ \  \ =\int_{U_xM}\tr_gs(\cdot,V,...,V,
\cdot,V,...,V)+...+\int_{U_xM}\tr_gs(V,...,V,\cdot, V,...,V, \cdot)
\end{array}
\end{equation}
for every $x\in M$.
\end{proposition}

We have the following  Ros' integral  formula
%\begin{equation}\label{RosI}
%\int_{UM}(\hat\nabla T)(v,...,v)=0,
%\end{equation}

\begin{equation}\label{RosII}
\int_{UM}\tr_g(\hat\nabla s)(\cdot,\cdot,V,...,V)=0,
\end{equation}
where $s$ is a covariant tensor field of degree greater than 1 on a
compact manifold $M$ and $\int_{UM}=\int_{x\in M}\int_{U_XM}$, see
\cite{R}.

\begin{thm} Let $(g,A)$ be a conjugate symmetric statistical
structure on a compact manifold $M$. If the sectional curvature
$\hat k\ge 0$ on $M$ and $\hat\nabla^2\tau =0$, then $\hat\nabla
A=0$ on $M$. If moreover, $\hat k>0$ at some point $p$ of $M$,  then
the statistical structure is trivial.
\end{thm}
\proof Consider the following tensor field $s$ on $M$
\begin{equation}
s(X_1,...,X_7)=\hat\nabla A(X_1,...,X_4)A(X_5,X_6,X_7).
\end{equation}
By Theorem \ref{nablasquaretau}, we have $\hat\nabla \tau=0$, which
implies that $\sum _i\hat\nabla A (X,Y,e_i,e_i)=0$, where, as usual,
$e_1,...,e_n$ stands for an orthonormal basis of  a tangent space.
Using (\ref{Ricci}), we obtain
\begin{equation}\label{pomocnicza1}
\begin{array}{rcl}
&&\sum_i\hat\nabla
s(e_i,e_i,V,...,V)\\
&&\ \ \ \ \ \ \
 \ \ \ \ \ \ =\sum_i\hat\nabla^2A(e_i,e_i,V,V,V)A(V,V,V)+\sum_i\hat\nabla
A(e_i,V,V,V)^2\\
&&\ \ \ \ \ \ \ \ \ \ \  \ \ =\Vert (\hat\nabla
K)(V,V,V)\Vert^2+\sum_i((\hat
R(e_i,V)A)(e_i,V,V)A(V,V,V)\\
&&\ \ \ \ \ \ =\Vert(\hat\nabla K)(V,V,V)\Vert ^2-\sum_iA(\hat
R(e_i,V)e_i,V,V)A(V,V,V)\\
&&\ \ \ \ \ \ \ \ \ \ \ \ \ \ \ \ \ \ \ \ \ \ \ \ \ \
-\sum_i2A(e_i,\hat R(e_i,V)V,V)A(V,V,V).
\end{array}
\end{equation}

If we  define a 1-form $\alpha$ on $U_pM$ by $ \alpha_V(e)=A(\hat
R(e,V)V,V,V)A(V,V,V) $, then, by (\ref{spherical_codifferential}),
we have
\begin{eqnarray*}
 &&\delta\alpha(V)=\sum_iA(\hat R(e_i,V)e_i,V,V)A(V,V,V)\\
&&\ \ \ \ \ \ \ \ \ \ +2\sum_iA(\hat R(e_i,V)V,e_i,V)A(V,V,V)\\
&&\ \ \ \ \ \ \ \ \ \ +3\sum _iA(\hat R(e_i,V)V,V,V)A(e_i,V,V),
\end{eqnarray*}
and consequently
\begin{eqnarray*}
&&\int_{U_pM}\left[-\sum_iA(\hat
R(e_i,V)e_i,V,V)A(V,V,V)-\sum_i 2A(e_i,\hat R(e_i,V)V,V)A(V,V,V)\right]\\
&&\ \ \ \ \ \ \ \ \  \ =3\int_{U_pM}g(\hat R (K(V,V),V)V,K(V,V)).
\end{eqnarray*}
By (\ref{RosII}) and (\ref{pomocnicza1}), one gets
\begin{equation}\label{ostatni}
0=\int_{UM}\Vert (\hat\nabla K)(V,V,V)\Vert ^2+3\int_{UM}g(\hat R
(K(V,V),V)V,K(V,V)).
\end{equation}
 This implies the first assertion. To
prove the second one, suppose that $\hat k>0$ at $p\in M$, $K\ne 0$
(automatically at every point, in particular at $p$). By
(\ref{ostatni}), $K(V,V)$ is parallel to $V$ for each $V\in U_pM$.
By
 Theorem \ref{nablasquaretau}, we know that $E$ vanishes on $M$.
Consider the function $U_pM\ni V\to A(V,V,V)\in \R$. It attains a
maximum at some vector, say $e_1$. It is an eigenvector of
$K_{e_1}$, that is, $K(e_1,e_1)=\lambda_1e_1$ and $\lambda_1>0$. Let
$e_1, ..., e_n$ be an orthonormal eigenbasis of $K_{e_1}$. Hence
$K(e_1, e_i)=\lambda _ie_i$ for some numbers $\lambda _i$. Since
$E=0$, we have $\sum _i\lambda _i=0$. Since $K(e_i,e_i)$ is parallel
to $e_i$ and $K(e_1,e_i)=\lambda_i e_i$, we have
$0=g(K(e_i,e_i),e_1)=g(K(e_1,e_i), e_i)=\lambda_i$  for every
$i=2,...,n$. It follows that $\lambda _1=0$, which gives a
contradiction.\koniec


\bigskip
\begin{thebibliography}{20}

\bibitem{Cal} Calabi E., \emph{Complete affine hypersurfaces I},
 Symposia Math. 10 (1972), 19-38.
\bibitem{ChBY46} Chen B.Y., Ogiue K., \emph{On totally rela submanifolds}
Trans. Amer. Math. Soc. 193 (1974), 257-266.

\bibitem{ChBY} Chen B.Y., \emph{ On Ricci curvature of isotropic and Lagrangian submanifolds in complex space forms,}
Arch. Math. (Basel) 74 (2000), 154-160.


\bibitem{CHBY1} Chen B.Y., \emph{Riemannian geometry of Lagrangian
submanifolds}, Taiwanese . J. Math., 5 (2001) 681-723.




\bibitem{DNV} Dillen F., Nomizu K., Vrancken L., \emph{Conjugate connections
 and Radon's theorem in affine differential geometry}, Monatsh.
 Math., 109 (1990) 221-235.


\bibitem{F} Furuhata H., Husegawa I., \emph{ Submanifold theory in
holomorphic statistical manifolds}, Geometry of Cauchy-Riemann
submanifolds, (2016), 179-215, Springer, Singapore.


\bibitem{K} Kurose T., \emph{Dual connections and affine geometry},
Math. Z.. 203 (1990), 115-121.

\bibitem{L} A-M. Li, \emph{Some theorems in Affine Differential
Geometry}, Acta M. Sci., 5 (1989), 345-354.

\bibitem{L_1} A-M. Li, \emph{ An intrinsic rigidity theorem for minimal
submanifolds in a sphere}, Arch. Math., 58 (1992), 582-594.

\bibitem{LSZ} A-M. Li, U. Simon, G. Zhao, \emph{Global Affine Differential Geometry of
Hypersurfaces}, Walter de Gruyter,  1993.


\bibitem{MRU} S. Montiel, A. Ros, F. Urbano,
\emph{ Curvature pinching and eigenvalue rigidity for minimal
 submanifolds},
 Math. Z., 191 (1986), 537-548.

\bibitem{NS} K. Nomizu, T. Sasaki, \emph{Affine Differential
Geometry}, Cambridge University Press, 1994.

\bibitem{BW4} Opozda B., \emph{Bochner's technique for statistical
structures}, Ann. Glob. Anal. Geom., 48 (2015) 357-395.

\bibitem{seccur}Opozda B., \emph{ A sectional curvature for statistical
structures},  Linear Algebra Appl.,497 (2016), 134-161.




\bibitem{calabi} Opozda B., \emph{Curvature bounded conjugate symmetric
statistical structures with complete
 metric}, Ann. Glob. Anal.
Geom., 55 (2019) 687-702.



\bibitem{R} A. Ros, \emph{ A characterization of seven compact
Kaehler sumbanifolds by holomorhic pinching}, Annales Math., 121
(1985), 377-382.


\bibitem{shima} Shima H. \emph{The geometry Hessian structure},
World Scientific Publishing, (2007).


\bibitem{Y} S.T. Yau, \emph{Harmonic functions on complete Riemannian
manifolds}, Comm. Pure Appl. Math., 28 (1975), 201-228.



\end{thebibliography}
\end{document}